\documentclass [12pt]  {article} 
\usepackage{amssymb}
\def \CC{{\mathbb{C}}}

\begin{document}

\begin{center}
{\Large {\bf What is the definition of two meromorphic
functions sharing a small function?}}\\
\bigskip
\bigskip
{\sc Andreas Schweizer}\\
\bigskip
{\small {\rm Department of Mathematics,\\
Korea Advanced Institute of Science and Technology (KAIST),\\ 
Daejeon 305-701\\
South Korea\\
e-mail: schweizer@kaist.ac.kr}}
\end{center}
\begin{abstract}
\noindent
Two meromorphic functions $f(z)$ and $g(z)$ sharing
a small function $\alpha(z)$ usually is defined in 
terms of vanishing of the functions $f-\alpha$ and 
$g-\alpha$. We argue that it would be better to modify 
this definition at the points where $\alpha$ has poles. 
Related to this issue we also point out some possible 
gaps in proofs in the published literature.
\\ 
{\bf Mathematics Subject Classification (2010):} 
30D35
\\
{\bf Key words:} meromorphic function; sharing a small function;
shared function; sharing IM; sharing CM
\end{abstract}

\noindent
The question in the title might surprise. After all, most 
papers on meromorphic functions that share a small function 
give a definition. Usually it is
\\ \\
{\bf Definition 1.} 
Let $f(z)$, $g(z)$ and $\alpha(z)$ be meromorphic functions 
on a domain $D$. We say that the functions $f$ and $g$ 
{\bf share the function} $\alpha$ {\bf in the sense of vanishing} 
if on $D$ the functions $f(z)-\alpha(z)$ and $g(z)-\alpha(z)$
have the same zeroes. More precisely we call this sharing the
function $\alpha$ {\bf IM (ignoring multiplicities) in the sense 
of vanishing}.
\par
If the zeroes of $f-\alpha$ and $g-\alpha$ not just coincide
in location but also in multiplicity, we say that $f$ and $g$ 
{\bf share $\alpha$ CM (counting multiplicities) in the sense 
of vanishing}.
\\ \\
So in short, $f$ and $g$ share $\alpha$ IM resp. CM in the
sense of vanishing if and only if the functions $f-\alpha$ 
and $g-\alpha$ share the value $0$ IM resp. CM.
\par
The words 'in the sense of vanishing' are not standard and 
were added by us to be able to distinguish this from another 
definition of sharing, which we will define below.
\par
Definition 1 ties in nicely with the generalizations of the
Second Main Theorem that involve the counting function of the 
zeroes of $f-\alpha$. 
See for example [7, Section 1.5, Theorem 1.36].
We refer to this book for all background information desired.
The above definition of sharing is also used in this book in
Sections 3.1.4 and 3.1.5.
\par
The smallness of $\alpha(z)$ is irrelevant for the definition 
of the sharing. In general the purpose of the smallness is to
allow Nevanlinna Theory to extract nontrivial consequences 
from the sharing.
\par
When we give examples in this paper, we always construct them 
in such a way that  $f$ and $g$ are meromorphic functions of 
finite order and $\alpha$ is a small function with respect 
to $f$ and $g$.
\\ \\
{\bf Example 1.}
The functions
$$f=\frac{1}{z}+e^z \ \ \ \hbox{\rm and}\ \  
g=\frac{1}{z}-\frac{e^z}{z}$$
share $\frac{1}{z}$ CM in the sense of vanishing. This is at
odds with intuition, because $\alpha$ and $f$ have a simple 
pole with residuum $1$ at $z=0$, whereas $g$ is entire.
\\ \\
We shall see that this is not the only problem.
\par
A useful fact is that if $f$ and $g$ share a value $a$ and 
if $M$ is a M\"obius transformation, then $M(f)$ and $M(g)$ 
share the value $M(a)$. In other words, sharing of values 
behaves well under M\"obius transformations. Obviously, 
translations and scaling do also respect sharing of a small 
function in the sense of vanishing. However, the inversion 
$z\mapsto \frac{1}{z}$ causes problems.
\\ \\
{\bf Example 2.}
The functions
$$f=z+z^2 e^z\ \ \ \hbox{\rm and}\ \  
g=z+z^3 e^z$$
share $\alpha=z$ IM in the sense of vanishing, 
but $\frac{1}{f}$ and $\frac{1}{g}$ do not share 
$\frac{1}{\alpha}=\frac{1}{z}$ in the sense of 
vanishing.
\par
However, the only point that causes trouble is $z=0$.
And since both, $\frac{1}{f}$ and $\frac{1}{g}$, have 
a simple pole with residuum $1$ at $z=0$, one sort of feels 
that they should share $\frac{1}{z}$ also at this point.
\\ \\
Even starting with two functions that share 
a small function CM does not help. 
\\ \\
{\bf Example 3.}
The functions
$$f=\frac{1}{z}+e^z\ \ \ \hbox{\rm and}\ \  
g=\frac{1}{z}-\frac{e^z}{z}$$
from Example 1 share $\alpha=\frac{1}{z}$ CM in the sense 
of vanishing, but $\frac{1}{f}$ and $\frac{1}{g}$ do not 
share $\frac{1}{\alpha}=z$ in the sense of vanishing, not 
even IM. Note that $\frac{1}{f}-z=\frac{-z^2 e^z}{1+ze^z}$
vanishes at $z=0$, but $\frac{1}{g}-z=\frac{ze^z}{1-e^z}$
takes the value $-1$ at $z=0$.
\\ \\
And even if the sharing survives the inversion, 
the multiplicities might not.
\\ \\
{\bf Example 4.}
The functions
$$f=\frac{1}{z}+e^z\ \ \ \hbox{\rm and}\ \  
g=\frac{1}{z}+\frac{e^z}{z}$$
share $\alpha=\frac{1}{z}$ CM in the sense of vanishing, 
but $\frac{1}{f}$ and $\frac{1}{g}$ share 
$\frac{1}{\alpha}=z$ only IM in the sense of vanishing.
\\ \\
We will come back to this. But let us first point out 
yet another problem.
\par
One argument that we have seen in several papers is the following.
Let $f$ and $g$ be two meromorphic functions and $\alpha$ a small 
function that is neither constant $0$ nor constant $\infty$. If
$f$ and $g$ share $\alpha$ CM (in the sense of vanishing), then 
$\frac{f}{\alpha}$ and $\frac{g}{\alpha}$ share the value $1$ CM. 
See for example the proofs of 
[2, Theorems 3,4, and 5], 
[3, Theorems 1.1 and 1.2],
[8, Theorems 1.1, 1.2, and 1.5], and
[9, Theorem 2].
\par
However, in general this claim is not true.
\\ \\
{\bf Example 5.}
Let 
$$f=\frac{1}{z}+e^z,\ \ \ 
g=\frac{1}{z}+\frac{e^z}{z},\ \ \ 
\hbox{\rm and}\ \  \alpha =\frac{1}{z}$$
as in Example 4.
Then $f$ and $g$ share $\alpha$ CM in the sense of vanishing. 
But $\frac{f}{\alpha}$ and $\frac{g}{\alpha}$ do {\bf not} 
share the value $1$, not even IM. (On $\CC^*$ they would share 
1 CM.)
\\ \\
One can even give examples where the sharing 
gets lost at infinitely many points.
\\ \\
{\bf Example 6.}
Let 
$$f=\frac{1}{\sin z}+e^{z^2},\ \ \ 
g=\frac{1+e^{z^2}}{\sin z},\ \ \ 
\hbox{\rm and}\ \  \alpha =\frac{1}{\sin z}.$$
Then $f$ and $g$ share $\alpha$ CM in the sense of vanishing, 
but $\frac{f}{\alpha}$ and $\frac{g}{\alpha}$ do {\bf not} 
share the value $1$, not even IM.
\\ \\
To overcome some of these problems, let us look at another 
possible (and in our opinion better) definition of sharing
a meromorphic function. 
\par
As is already evoked by the expression 'moving target', one
can think of $\alpha(z)$ as a value that is changing with $z$. 
\\ \\
{\bf Definition 2.} 
Let $f(z)$, $g(z)$ and $\alpha(z)$ be meromorphic functions 
on a domain $D$. We say that the functions $f$ and $g$ 
{\bf share the function} $\alpha$ 
{\bf (IM) in the sense of value} if for every $z_0 \in D$ 
we have
$$f(z_0)=\alpha(z_0)\Leftrightarrow g(z_0)=\alpha(z_0).$$
At the points $z_0$ where $\alpha$ has a pole, 
$f(z_0)=\alpha(z_0)$ simply means that $f$ also has a pole 
at $z_0$.
\par
If $\alpha$ is constant (including the possibility 
$\alpha\equiv\infty$) this definition specializes to
the usual definition of sharing a value IM.
\par
By the new definition, the functions $f$ and $g$ from 
Example 1 do indeed not share the function $\frac{1}{z}$.
\par
In the first version of this paper (arXiv:1705.05048v1) 
we had also given an ad-hoc definition of sharing a function
$\alpha$ CM in the sense of value. But that definition had 
several drawbacks; so we omit it here.
\par
Instead, we now give what we believe is the good definition 
of sharing a function CM. We will back up this conviction
by showing that this definition has some desirable properties 
that Definition 1 is lacking.
\par
We start from the observation that if $\alpha$ has a pole at 
$z_0$, then the condition $f(z_0)=\alpha(z_0)$ is equivalent
to the vanishing of $\frac{1}{f}-\frac{1}{\alpha}$ at $z_0$. 
This is the only modification made to get from Definition 
1 to Definition 2. Refining this, towards sharing CM, we 
take into account multiplicities, namely: outside the poles 
of $\alpha$ the order of vanishing of $f-\alpha$ and $g-\alpha$
as before, and at the poles of $\alpha$ the order of vanishing
of $\frac{1}{f}-\frac{1}{\alpha}$ and 
$\frac{1}{g}-\frac{1}{\alpha}$.
\\ \\
{\bf Definition 3.} 
Let $f(z)$, $g(z)$ and $\alpha(z)$ be meromorphic functions 
on a domain $D$. We say that the functions $f$ and $g$ 
{\bf share the function} $\alpha$ {\bf CM in the sense of value} 
if for every $z_0 \in D$ the following conditions hold:
\begin{itemize}
\item
At points $z_0$ with $\alpha(z_0)\neq\infty$: The function 
$f-\alpha$ has a zero of order $m$ at $z_0$ if and only if
$g-\alpha$ has a zero of order $m$ at $z_0$.
\item
At points $z_0$ with $\alpha(z_0)=\infty$: The function 
$\frac{1}{f}-\frac{1}{\alpha}$ has a zero of order $m$ at 
$z_0$ if and only if $\frac{1}{g}-\frac{1}{\alpha}$ has 
a zero of order $m$ at $z_0$.
\end{itemize}
\bigskip
\noindent
Quite likely this definition has already been discussed and 
used in the literature. But we couldn't find any sources. And
an exhaustive search is hardly feasible given the sheer mass
of publications in this area.
\par
Comparing sharing in the sense of vanishing and sharing in the 
sense of value, one obviously expects some disagreement at the
points where $\alpha(z)$ has a pole.
\par
We already mentioned Example 1, where $f$ and $g$ share 
$\frac{1}{z}$ in the sense of vanishing but not in the sense 
of value. The converse also occurs.
\\ \\
{\bf Example 7.}
The functions
$$f=\frac{1}{z}+ e^z \ \ \ \hbox{\rm and}\ \  
g=\frac{1}{z}+ z e^z $$
share $\frac{1}{z}$ IM in the sense of value (not CM!).
Note that $\frac{1}{f}-\frac{1}{\alpha}=\frac{-z^2 e^z}{1+ze^z}$
and $\frac{1}{g}-\frac{1}{\alpha}=\frac{-z^3 e^z}{1+z^2 e^z}$.
But our main point is that $f$ and $g$ do not share $\frac{1}{z}$ 
in the sense of vanishing.
\\ \\
Even if $f$ and $g$ share $\alpha$ in the sense of vanishing 
and in the sense of value, it is not guaranteed that the notions 
of sharing CM agree.
\\ \\
{\bf Example 8.}
Let 
$$f=\frac{1}{z}+e^z,\ \ \ 
g=\frac{1}{z}+\frac{e^z}{z},\ \ \ 
\hbox{\rm and}\ \  \alpha =\frac{1}{z}.$$
Then $f$ and $g$ share $\alpha$ CM in the sense of vanishing,
but only IM in the sense of value. Compare Example 4. The 
discrepancy occurs at $z=0$.
\\ \\
Here are two reasons why sharing in the sense of value is better.
\\ \\
{\bf Theorem 1.} \it 
Sharing in the sense of value (IM or CM) is well-behaved under
M\"obius transformations. More precisely:
\par
If $f$ and $g$ share $\alpha$ (IM or CM) in the sense of value
and if $M$ is a M\"obius transformation, then $M(f)$ and $M(g)$
share $M(\alpha)$ (IM resp. CM) in the sense of value.
\rm
\\ \\
{\bf Proof.}
Obviously, translations and scaling respect any type of sharing 
we have discussed so far. So it suffices to prove the statement
of the theorem for the inversion $z\mapsto\frac{1}{z}$.
\par
Let us first look at points $z_0$ with $\alpha(z_0)\in\CC^*$.
From $\frac{1}{f}-\frac{1}{\alpha}=\frac{\alpha-f}{\alpha f}$
we then see that $f-\alpha$ has a zero of order $m$ ($>0$) at
$z_0$ if and only if $\frac{1}{f}-\frac{1}{\alpha}$ has a zero 
of order $m$ at $z_0$. Since the same holds of course for $g$,
we obtain that $z\mapsto\frac{1}{z}$ respects the sharing 
(IM or CM) outside the zeroes and poles of $\alpha$. 
\par
At the points with $\alpha(z_0)=0$ or $\alpha(z_0)=\infty$
this also holds by definition of sharing in the sense of value.
\hfill $\Box$
\\ \\
{\bf Theorem 2.} \it
Let $f$, $g$ and $\alpha$ be meromorphic functions on a domain 
$D$, such that $\alpha$ is not constant $0$ and also not constant 
$\infty$. Suppose that on $D$ the functions $f$ and $g$ share 
the function $\alpha$ CM in the sense of value. Then 
$\frac{f}{\alpha}$ and $\frac{g}{\alpha}$ share the value 
$1$ CM on $D$.
\rm
\\ \\
{\bf Proof.}
At points $z_0$ with $\alpha(z_0)\in\CC^*$ the functions $f-\alpha$ 
and $\frac{f}{\alpha}-1$ obviously have the same order of vanishing. 
\par
Now let $\alpha(z_0)=0$. By assumption, $f-\alpha$ and $g-\alpha$
either both have a zero of the same order at $z_0$ or both do not
vanish at $z_0$. So the functions 
$\frac{f-\alpha}{\alpha}=\frac{f}{\alpha}-1$ and
$\frac{g-\alpha}{\alpha}=\frac{g}{\alpha}-1$ either both have 
a zero of the same order at $z_0$ or both do not vanish at $z_0$.
This is what we wanted.
\par
At points $z_0$ with $\alpha(z_0)=\infty$ we use that by
Theorem 1 the functions $\frac{1}{f}$ and $\frac{1}{g}$ share
$\frac{1}{\alpha}$ CM in the sense of value. Then by the previous
argument $\frac{\alpha}{f}$ and $\frac{\alpha}{g}$ share the
value $1$ CM.  Another application of $z\mapsto\frac{1}{z}$ 
finishes the proof.
\hfill $\Box$
\\ \\
However, this does not mean that one can get rid of all problems 
concerning $\frac{f}{\alpha}$ and $\frac{g}{\alpha}$ sharing
$1$ CM by simply changing the definition of sharing a small 
function from sharing in the sense of vanishing to sharing in the 
sense of value. This is a change that comes at a certain price
and could create new problems somewhere else.
For example, if $f$ and $g$ share $\alpha$ in the sense of 
vanishing, every zero of $f-\alpha$ also is a zero of $g-\alpha$
and of $f-g$. But from sharing in the sense of value we only get
that every zero of $f-\alpha$ that is not a pole of $\alpha$ also
is a zero of $g-\alpha$ and of $f-g$.
\par
Another claim that can be found in the literature is that
if two meromorphic functions $f$ and $g$ share the small 
function $\alpha$ ($\not\equiv 0$, $\not\equiv\infty$) IM, 
then $\frac{f}{\alpha}$ and $\frac{g}{\alpha}$ share the 
value $1$ IM. See for example the proofs of 
[3, Theorems 1.1 and 1.2],
[5, Theorems 1 and 2], and
[9, Theorem 1].
\par
Again the claim is not true in general. Just see 
Example 5 ($=$ Example 8), where $f$ and $g$ share $\alpha$ 
in the sense of vanishing and in the sense of value, but 
$\frac{f}{\alpha}$ and $\frac{g}{\alpha}$ do not share the 
value $1$.
\par
But one can give a more striking counter-example in which 
all functions involved are entire.
\\ \\
{\bf Example 9.}
Let 
$$f=(\sin z) +(\sin z) e^{z^2},\ \ \ 
g=(\sin z) +(\sin z)^2 e^{z^2},\ \ \ 
\hbox{\rm and} \ \ \alpha =\sin z.$$
Then $f$ and $g$ share $\alpha$ IM in the sense of vanishing 
{\bf and} in the sense of value, but $\frac{f}{\alpha}$ 
and $\frac{g}{\alpha}$ do {\bf not} share the value $1$.
\\ \\
In this case the problems are coming from the zeroes of 
$\alpha$. This example also shows that there can hardly 
be a reasonable way to define the problem away.
\par
We do of course not claim that our Examples 5, 6 and 
9 are direct counterexamples to the specified theorems 
in the papers we mentioned. 
The functions in these papers satisfy many more conditions, 
for example $g$ being a derivative of $f$ or a differential 
polynomial of $f$, or there are additional sharing properties.
But without additional arguments the status of those theorems
is questionable.
\par
Since $\alpha$ is small, the counting function of the poles 
and zeroes of $\alpha$ is small, even when they are counted 
with their multiplicities. (This is not necessarily true when
they are counted with their multiplicities as poles or zeroes 
of $f$, as those multiplicities might grow rapidly.) So if 
$f$ and $g$ share $\alpha$ IM (resp. CM), one can still say
that the truncated counting function of the points where
$\frac{f}{\alpha}$ and $\frac{g}{\alpha}$ do not share the value
$1$ IM (resp. not share the value  $1$ CM) is small. 
\par
But it makes a difference whether one can apply a well-known
theorem on functions that share the value $1$ CM (compare for
example [4, Theorem A]), or whether one would need a theorem 
on functions that share the value $1$ outside a small set of
arguments.
\par
Many papers argue correctly that $\frac{f}{\alpha}$ and 
$\frac{g}{\alpha}$ share the value $1$ outside the zeroes 
and poles of $\alpha$.
Some other papers avoid the whole problem by imposing the 
extra condition that $\alpha(z)$ has no common poles and 
no common zeroes with $f(z)$ or $g(z)$. 
\par
Before we continue our discussion of this, we have to recall
the definition of weighted sharing, which has proved to be very 
useful in the theory of value sharing. Weighted sharing (of
values) was introduced by Lahiri in [4] to have some finer 
degrees of division between sharing CM and sharing IM.
\\ \\
{\bf Definition 4.} 
Let $f(z)$, $g(z)$ and $\alpha(z)$ be meromorphic functions 
on a domain $D$. We say that the functions $f$ and $g$ 
{\bf share the function} $\alpha$ {\bf with weight $m$ 
in the sense of vanishing} if for each $k=1,2,\ldots,m$
the $k$-fold zeroes of $f-\alpha$ coincide with the 
$k$-fold zeroes of $g-\alpha$, and the zeroes of $f-\alpha$
of multiplicity bigger than $m$ coincide with the zeroes of
$g-\alpha$ of multiplicity bigger than $m$. 
\\ \\
In the latter case the multiplicities are not necessarily 
the same. So sharing with weight $\infty$ is sharing CM, 
and sharing with weight $0$ is sharing IM.
\par
Analoguously we could refine Definition 3 into a definition 
of sharing $\alpha$ with weight $m$ in the sense of value.
By Theorem 1 and its proof at points $z_0$ with 
$\alpha(z_0)\in\CC^*$ a $k$-fold zero of $f-\alpha$ at $z_0$
corresponds to a $k$-fold zero of $\frac{1}{f}-\frac{1}{\alpha}$
at $z_0$. So we can build this into the formulation and give
the following equivalent definition.
\\ \\
{\bf Definition 5.} 
Let $f(z)$, $g(z)$ and $\alpha(z)$ be meromorphic functions 
on a domain $D$. We say that the functions $f$ and $g$ 
{\bf share the function} $\alpha$ {\bf with weight $m$ 
in the sense of value} if outside the poles of $\alpha$ 
the functions $f$ and $g$ share $\alpha$ with weight $m$ 
in the sense of vanishing, and outside the zeroes of 
$\alpha$ the functions $\frac{1}{f}$ and $\frac{1}{g}$ 
share $\frac{1}{\alpha}$ with weight $m$ in the sense 
of vanishing.
\\ \\
This definition contains everything we need as special cases,
namely sharing IM in the sense of value (weight $0$), sharing
CM in the sense of value (weight $\infty$), and sharing a value
with weight. Also, when $\alpha$ has no poles, sharing $\alpha$
(IM, CM, with weight $m$) in the sense of value coincides with
the corresponding notion of sharing $\alpha$ in the sense of 
vanishing.
\par
As we already pointed out, Definition 5 is well-behaved under 
M\"obius transformations.
\par
The functions $f$ and $g$ in Example 7 share $\frac{1}{z}$ 
with weight $1$ in the sense of value.
\par
The proof of [6, Theorem 1.1] contains a statement of the form
that if $f$ and $g$ share $\alpha$ with weight $m$, then
$\frac{f}{\alpha}$ and $\frac{g}{\alpha}$ share the value
$1$ with weight $m$. 
\par
Unfortunately, this claim is also not true in general,
not even if all functions involved are entire. To see
that we generalize Example 9.
\\ \\
{\bf Example 10.}
Let $m$ be a nonnegative integer and
$$f=(\sin z)^{m+1} +(\sin z)^{m+1} e^{z^2},\ \ \ 
g=(\sin z)^{m+1} +(\sin z)^{m+2} e^{z^2},\ \ \ 
\hbox{\rm and}\ \  \alpha =(\sin z)^{m+1}.$$
Then $f$ and $g$ share $\alpha$ with weight $m$ 
in the sense of vanishing {\bf and} in the sense 
of value, but $\frac{f}{\alpha}$ and $\frac{g}{\alpha}$ 
do {\bf not} share the value $1$, not even IM.
\\ \\
In this case the problems are coming from the zeroes of 
$\alpha$. And again we see that there is probably no hope
of completely getting rid of the problem.
\par
In contrast, the claim in the proof of [1, Theorem 1.1]
that if $f$ and $g$ weakly share $\alpha$ with weight $2$,
then $\frac{f}{\alpha}$ and $\frac{g}{\alpha}$ share the value
$1$ weakly with weight $2$ is correct, because it is already
built into the definition of sharing weakly [1, Definition 1.3]
that there might be a small set where the sharing doesn't hold.
\par
Finally we point out that the converse of the claims discussed
above is also not true, not even if most of the functions 
involved are entire.
\\ \\
{\bf Example 11.}
Let
$$f=1 + e^{z^2}\ \ \ \hbox{\rm and}\ \  
g=1+\frac{e^{z^2}}{\sin z}.$$ 
Then $f$ and $g$ share the value $1$ CM, 
but $\alpha f$ and $\alpha g$ do {\bf not} share the small 
function $\alpha=\sin z$, neither in the sense of vanishing 
nor in the sense of value.
\\

\subsection*{\hspace*{10.5em} References}
\begin{itemize}

\item[{[1]}] A.~Banerjee and S.~Mukherjee: \rm Nonlinear
differential polynomials sharing a small function, \it
Archivum Mathematicum (Brno) \bf 44 \rm (2008), 41-56 

\item[{[2]}] K.~Boussaf, A.~Escassut, J.~Ojeda: \rm 
Complex meromorphic functions $f'P'(f)$, $g'P(g)$ sharing
a small function, \it Indag. Math. \bf 24 \rm (2013), 15-41 

\item[{[3]}] A.~Chen, X.~Wang, G.~Zhang: \rm Unicity of 
meromorphic function sharing one small function with its derivative, 
\it Int. J. Math. Math. Sci. \rm 2010 Art.ID 507454, 11pp.

\item[{[4]}] I.~Lahiri: \rm Weighted value sharing and uniqueness
of meromorphic functions, \it Complex Variables \bf 46 \rm
(2001), 241-253 

\item[{[5]}] S.~Majumder: \rm Unicity of meromorphic function
that share a small function with its derivative, \it 
Appl. Math. E-Notes \bf 14 \rm (2014), 144-150 

\item[{[6]}] P.~Sahoo and B.~Saha: \rm Uniqueness of meromorphic
functions whose certain differential polynomials share a small
function, \it Acta Math. Vietnam \bf 41 \rm (2016), 25-36

\item[{[7]}] C.-C.~Yang and H.-X.~Yi: \it Uniqueness theory
of meromorphic functions, \rm Kluwer Academic Publishers Group,
Dordrecht, 2003

\item[{[8]}] X.-B.~Zhang and J.-F.~Xu: \rm Uniqueness of 
meromorphic functions sharing a small function and its applications,
\it Comput. Math. Appl. \bf 61 \rm (2011), 722-730 

\item[{[9]}] J.-L.~Zhang and L.-Z.~Yang: \rm Some results 
related to a conjecture of R.~Br\"uck concerning meromorphic
functions sharing one small function with their derivatives,
\it Ann. Acad. Sci. Fenn. Math. \bf 32 \rm (2007), 141-149 

\end{itemize}

\end{document}